%% file: JVR.tex
\documentclass[12pt]{article}

\usepackage{epsfig}
\usepackage{amsopn,amsthm,amssymb,amsmath}

\theoremstyle{plain}
\newtheorem{Thm}{Theorem}[section]
\newtheorem{Prop}[Thm]{Proposition}
\newtheorem{Lem}[Thm]{Lemma}
\newtheorem{Cor}[Thm]{Corollary}

\theoremstyle{definition}
\newtheorem*{Def}{Definition}

\theoremstyle{remark}

\newtheorem{Remark}[Thm]{Remark}

\newcommand{\nn}{\mathbb{N}}
\newcommand{\zz}{\mathbb{Z}}
\newcommand{\rr}{\mathbb{R}}
\newcommand{\e}{\varepsilon}
\newcommand{\ph}{\varphi}
\newcommand{\oo}{\mathcal{O}}
\newcommand{\nnn}{\mathcal{N}}
\newcommand{\rmap}{\longrightarrow}
\newcommand{\rrt}{\widetilde{\rr}}
\newcommand{\dvec}[1]{\displaystyle\overrightarrow{#1}}
\newcommand{\ddvec}[1]{\overrightarrow{#1}}
\newcommand{\xxangle}{\theta}
\newcommand{\xxTN}{\mathcal{P}}

\newcommand{\xxDiv}{\operatorname{Div}}
\newcommand{\xxJac}{\operatorname{Jac}}
\newcommand{\xxBunch}{\operatorname{Bunch}}
\newcommand{\xxInt}{\operatorname{Int}}

\newcommand{\xxwt}{\operatorname{wt}}
\newcommand{\xxNewt}{\operatorname{Newt}}
\newcommand{\xxConv}{\operatorname{Conv}}
\newcommand{\xxVer}{\operatorname{Ver}}
\newcommand{\xxEdge}{\operatorname{Edge}}
\newcommand{\xxmoment}{\operatorname{moment}}
\newcommand{\xxMult}{\operatorname{Mult}}
\newcommand{\xxarea}{\operatorname{area}}
\newcommand{\xxlength}{\operatorname{length}}
\newcommand{\xxalg}{\operatorname{alg}}

\title {Jacobian varieties of reduced tropical curves}
\author {Shuhei Yoshitomi}
\date{}
\begin {document}
\maketitle

\begin {abstract}

On tropical geometry in $\rr^2$, the divisor and the Jacobian variety are defined in analogy to algebraic geometry. For study of these objects, it is important to think of the `bunch' of a tropical curve (Figure 1). In this paper, we will show that if the bunch is a bouquet, then the Jacobian is a higher-dimensional torus.

\end {abstract}

\begin {figure}[ht]
\begin {center} \input {./picture/bhA1.tex} \end {center}
\caption {Bunch of a tropical curve}
\label {Fig:A1}
\end {figure}
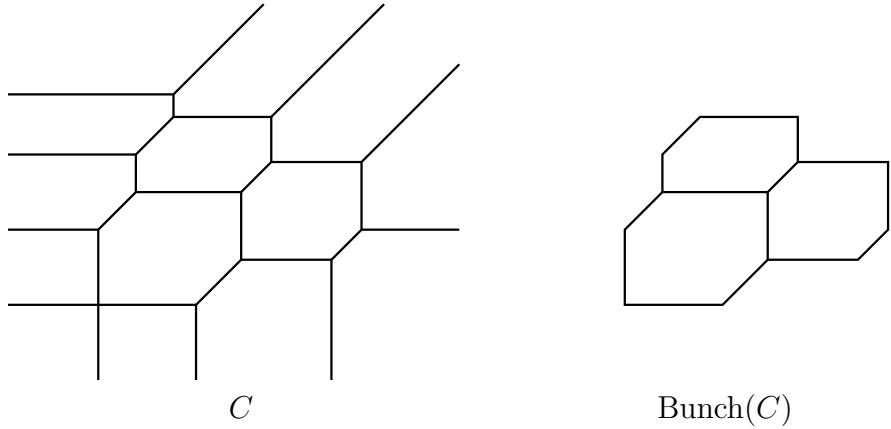

\section {Introduction}

In this paper, the affine space $\rr^2$ is equipped with interior product $(u_1, u_2) \cdot (v_1, v_2)= u_1 v_1+ u_2 v_2$ and exterior product $(u_1, u_2) \times (v_1, v_2)= u_1 v_2- v_1 u_2$. \par
Let $C$ be a tropical curve in $\rr^2$ (See section \ref{Sec:B1} for preliminary). The Jacobian variety of $C$ is defined in analogy to algebraic geometry as follows.

\begin {Def}
The {\em divisor group} $\xxDiv(C)$ of $C$ is the free abelian group generated by all points of $C$. We define a subgroup
\[
\xxDiv^0(C) = \left\{ D= \sum_{P \in C}^{} m_P P \in \xxDiv(C) ~|~ \deg D := \sum_P m_P= 0 \right\}.
\]
Divisors $D, D'$ are {\em linearly equivalent}, $D \sim D'$, if there are tropical curves $L, L'$ such that
\[
\Delta(L)= \Delta(L'), \]\[
D-D'= C \cdot L- C \cdot L',
\]
where $\Delta$ denotes the Newton polygon and $C \cdot L$ denotes the stable intersection. \\
The {\em Jacobian variety} of $C$ is the residue group
\[
\xxJac(C)= \xxDiv^0(C) / \sim.
\]
\end {Def}

\begin {figure}[ht]
\begin {center} \input {./picture/bhA3.tex} \end {center}
\caption {Tropical elliptic curve}
\label {Fig:A3}
\end {figure}
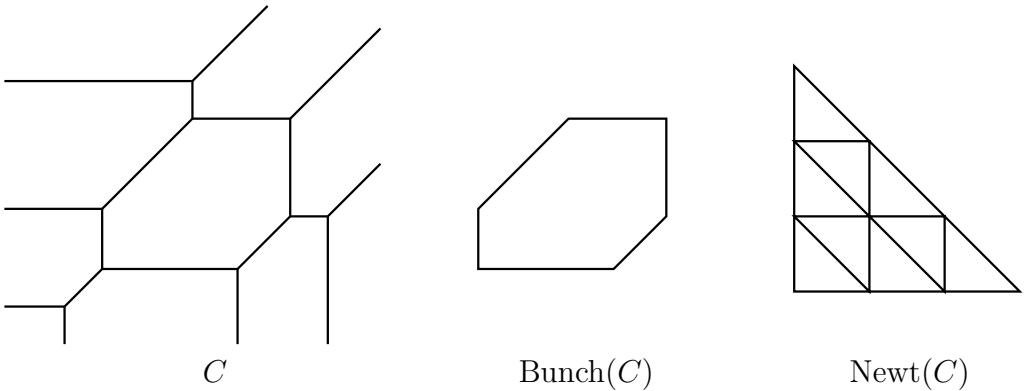

We fix a point $\oo \in C$. If $C$ is a tropical elliptic curve, then the `bunch' of $C$ is homeomorphic to $S^1$ (Figure \ref{Fig:A3}), and we have a map
\begin {eqnarray*}
\ph \colon \xxBunch(C) & \rmap & \xxJac(C) \\
P & \mapsto & P- \oo.
\end {eqnarray*}
Vigeland \cite{elliptic curve} states that $\ph$ is bijective, which concludes that a tropical elliptic curve has a group structure in a smaller part of it. But the proof in \cite{elliptic curve} is not complete. He has proved the surjectivity of $\ph$. \par
In this paper, we give a complete proof and some generalization of it.

\begin {Def}
A tropical curve $C$ is {\em reduced} if every edge is of weight $1$. $C$ is {\em smooth} if every vertex is $3$-varent and of multiplicity $1$. (Hence any smooth tropical curve is reduced.)
\end {Def}

\begin {Def}
An edge $E \subset C$ is {\em tentacle} if $C \setminus \xxInt(E)$ is disconnected. $E$ is a {\em ray} if $E$ has only one vertex. We define the {\em bunch} of $C$, $\xxBunch(C)$, to be the quotient space of $C$ by every tentacle edge and ray contracted.
\end {Def}

\begin {Def}
A {\em bouquet} is a topological space $B= \Lambda_1 \cup \cdots \cup \Lambda_g$ with a point $\oo \in B$ such that
\[
\Lambda_i \approx S^1 \quad (1 \leq i \leq g), \]\[
\Lambda_i \cap \Lambda_j= \{ \oo \} \quad (i \not= j).
\]
$\oo$ is called the {\em center} of $B$. (The topological genus of $B$ is $g$.)
\end {Def}

\begin {Thm} \label {Prop:A1}
Let $C$ be a reduced tropical curve in $\rr^2$. Suppose that $\xxBunch(C)$ is a bouquet of genus $g$, with center $\oo$ and cycles $\Lambda_1, \ldots, \Lambda_g$. Then the map
\begin {eqnarray*}
\ph \colon \Lambda_1 \times \cdots \times \Lambda_g & \rmap & \xxJac(C) \\
(P_1, \ldots, P_g) & \mapsto & P_1+ \cdots+ P_g- g \oo
\end {eqnarray*}
is bijective.
\end {Thm}

The statement includes that the map $\ph$ is well-defined, i.e. if $P_1, P_1'$ are on the same tentacle edge or on the same ray, then $P_1 \sim P_1'$.

\begin {Remark}
Tropical geometry is introduced in three approaches. The first is an algebraic approach (e.g.\cite{Izh},\cite{Stu}). If $f \in \rrt[x, y]$ is a tropical polynomial over the tropical algebra $\rrt = \rr \cup \{ - \infty \}$, then its corner locus is a tropical curve. The second is a valuation-theoretical approach (e.g.\cite{Gath},\cite{Mik}). If $K$ is a suitable valuation field, and $V$ is an algebraic curve in $K^2$, then its image by the valuation map $v \colon K^2 \rightarrow \rrt^2$ is a tropical curve. This image is expressed as the limit of amoebas of complex curves. And the last is a geometrical approach (e.g.\cite{Mik}). A $1$-dimensional simplicial complex in $\rr^2$ (or, a graph in $\rr^2$) satisfying some balancing condition is a tropical curve. The equivalence of these definitions is easy to prove (See \cite{Gath}). In this paper, we take only geometrical approach. \par
\end {Remark}

\begin {Remark}
In algebraic approach, projective tropical curves are treated as special objects ($C$ is {\em projective} if $\Delta(C)= \Delta_d$). The above definition of linearly equivalence has another version as follows: $D \sim_{\xxalg} D'$ if there are projective tropical curves $L, L'$ such that
\[
\deg(L)= \deg(L'), \]\[
D-D'= C \cdot L- C \cdot L'.
\]
The residue group $\xxalg\xxJac(C)= \xxDiv^0(C) / \sim_{\xxalg}$ will be called the {\em algebraic Jacobian variety} of $C$. There is a canonical surjection $\psi \colon \xxalg\xxJac(C) \rightarrow \xxJac(C)$. We can think of a map $\ph_{\xxalg} \colon \xxBunch(C) \rightarrow \xxalg\xxJac(C)$ instead of $\ph \colon \xxBunch(C) \rightarrow \xxJac(C)$ (in the case of genus $1$). Vigeland \cite{elliptic curve} exactly states that $\ph_{\xxalg}$ is bijective. The injectivity of $\ph_{\xxalg}$ follows from the injectivity of $\ph$, but the surjectivity of $\ph_{\xxalg}$ requires extra arguments like \cite{elliptic curve}.
\end {Remark}

\begin {figure}[ht]
\begin {center} \input {./picture/bhA4.tex} \end {center}
\caption {}
\label {Fig:A4}
\end {figure}

\begin {Remark}
One needs to be careful about the `degree' of the tropical curve. Tropical curves $L, M$ with Newton polygons
\[
\Delta(L)= \xxConv \{(0, 0), (1, 0), (1, 1) \}, \]\[
\Delta(M)= \xxConv \{(0, 0), (0, 1), (1, 1) \}
\]
are both `tropical curves of degree $2$'. Figure \ref{Fig:A4} is an example such that $C \cdot L- C \cdot M= P- Q$. But Theorem \ref{Prop:A1} asserts that $P, Q$ are not linearly equivalent. Note that $L, M$ are not projective.
\end {Remark}

\section {Preliminary} \label {Sec:B1}

A {\em primitive vector} in $\rr^2$ is an integral vector $u= (u_1, u_2)$ such that $u_1, u_2$ are relatively prime. Any integral vector $v \in \zz^2$ is a primitive vector $u$ times some natural number $m$. $m$ is called the {\em lattice length} of $v$. \par
Let $C$ be a $1$-dimensional simplicial complex of rational slopes in $\rr^2$. Each finite edge $E \subset C$ has two primitive vectors for the directions parallel to $E$, say $u, -u$. If the {\em weight} $m \in \nn$ of $E$ is given, we have the weighted primitive vectors $mu, -mu$ of $E$. For a vertex $V \in E$, one of these vectors has the direction from $V$ to $E$. We call this, say $u_E$, the weighted primitive vector of $E$ starting at $V$. \par
If $E$ is a ray (i.e. an infinite edge), $E$ has only one weighted primitive vector.

\begin {Def}
$V$ {\em satisfies the balancing condition} if the sum of all weighted primitive vectors starting at $V$ equals $0$:
\[
\sum_{E \ni V}^{}u_E = 0.
\]
A {\em tropical curve} in $\rr^2$ is a $1$-dimensional weighted simplicial complex of rational slopes such that each vertex satisfies the balancing condition.
\end {Def}

Let $C$ be a tropical curve in $\rr^2$. Let $U_1, \ldots , U_r$ be all connected components of $\rr^2 \setminus C$. Let $N$ be a $1$-dimensional simplicial complex with vertex set $\xxVer(N)= \{w_1, \ldots, w_r\}$, $w_i \in \zz^2$.

\begin {Def}
$N$ is a {\em Newton complex} of $C$ if it satisfies the following conditions for any $i \not= j$. \\
i) If $\overline{U_i} \cap \overline{U_j}= \emptyset$, then $[w_i, w_j] \notin \xxEdge(N)$. \\
ii) If $\overline{U_i} \cap \overline{U_j}= E$ for some $E \in \xxEdge(C)$, then $[w_i, w_j] \in \xxEdge(N)$, and $w_j- w_i$ has lattice length $\xxwt(E)$, direction orthogonal to $E$ from $U_i$ to $U_j$.
\end {Def}

\begin {Prop} \label {Thm:A1}
Let $C, C_1, C_2$ be tropical curves in $\rr^2$. \\
1) A Newton complex $\xxNewt(C)$ of $C$ exists uniquely up to parallel translation. The convex hull
\[
\Delta(C)= \xxConv( \xxNewt(C))
\]
is called the {\em Newton polygon} of $C$. \\
2) $\Delta(C_1 \cup C_2)$ equals the Minkowski sum of $\Delta(C_1)$ and $\Delta(C_2)$.
\end {Prop}

\begin {figure}[ht]
\begin {center} \input {./picture/bhA2.tex} \end {center}
\caption {Newton complex}
\label {Fig:A2}
\end {figure}
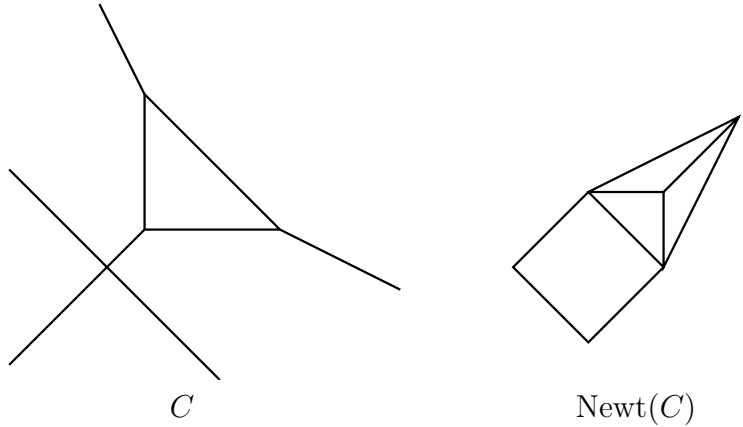

\begin {proof}
1) (See \cite{Mik}, \S3.4, if you take an algebraic approach.) Let $w_1= (0, 0)$, and suppose $w_1, \ldots, w_{k-1}$ are constructed. Rearranging $U_k, \ldots, U_r$, we may assume $U_i \cap U_k= E_k$ for some $i< k$ and some $E_k \in \xxEdge(C)$. Let $u_k$ be the primitive vector of direction orthogonal to $E_k$ from $U_i$ to $U_k$. Let $w_k$ be the vector
\begin {equation}
w_k- w_i= \xxwt(E_k) u_k. \label {eq:T-A1-1}
\end {equation}
If $U_1, \ldots, U_k$ are adjacent at a common vertex $V$, the condition (\ref{eq:T-A1-1}) is compatible for $U_k$ and $U_1$, i.e.
\[
w_1- w_k= \xxwt(E_1) u_1,
\]
where $E_1$ is the boundary of $U_k$ and $U_1$, and $u_1$ is the primitive vector of direction orthogonal to $E_1$ from $U_k$ to $U_1$. This compatibility follows from the balancing condition
\[
\sum_{i=1}^{k} \xxwt(E_i) u_i= 0.
\]
Therefore this construction does not depend on the choice of $U_k$. \par
2) $\Delta(C_1)$ depends only on the data of rays of $C_1$. Rays of $C_1 \cup C_2$ corresponds to rays of $C_1$ and $C_2$.
\end {proof}

$\xxNewt(C)$ can be considered as a dual object of $C$, with correspondence from $U_i$ to $w_i$ (Figure \ref{Fig:A2}). A vertex $V \in C$ corresponds to a polygon $T_V \subset \Delta(C)$ as follows. \\
i) $U_{i_1}, \ldots, U_{i_k}$ are adjacent at $V$, \\
ii) $T_V= \xxConv \{w_{i_1}, \ldots, w_{i_k} \}$.

\begin {Prop}[Global balancing condition] \label {Thm:A2}
Let $C$ be a tropical curve in $\rr^2$. Let $\Lambda$ be a simple closed curve in $\rr^2$ intersecting edges of $C$, say $E_1, \ldots, E_N$, transversely. Then
\[
\sum_{i=1}^{N} u_{E_i}= 0,
\]
where $u_{E_i}$ is the weighted primitive vector of $E_i$ starting at the vertex inside $\Lambda$.
\end {Prop}

\begin {proof}
For each vertex $V_j \in C$ inside $\Lambda$, the balancing condition
\[
\sum_{k} u_{jk}= 0
\]
holds. Thus
\[
\sum_{j,k} u_{jk} = 0.
\]
On the left hand side, two weighted primitive vectors of the same edge inside $\Lambda$ are canceled. Thus we have the required equation.
\end {proof}

A {\em tangent vector} $(v, P)$ in $\rr^2$ is a vector $v \in \rr^2$ with a starting point $P \in \rr^2$. We fix a point $P_0 \in \rr^2$. The {\em moment} of $(v, P)$ is the exterior product
\[
\xxmoment(v, P)= \ddvec{P_0 P} \times v.
\]

\begin {figure}[ht]
\begin {center} \input {./picture/bhB1.tex} \end {center}
\caption {Moment condition inside $\Lambda$. Note that $\xxmoment(u_{E_1}, V_{E_1}) \not= \xxmoment(u_{E'_1}, V_{E'_1})$.}
\label {Fig:B1}
\end {figure}
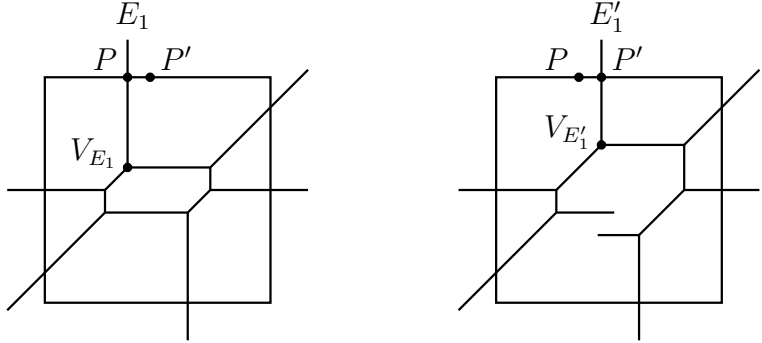

\begin {Prop}[Moment condition] \label {Thm:A3}
Under the assumption of Proposition \ref{Thm:A2},
\[
\sum_{i=1}^{N} \xxmoment(u_{E_i}, V_{E_i})= 0,
\]
where $(u_{E_i}, V_{E_i})$ is the weighted primitive tangent vector of $E_i$ starting at the vertex inside $\Lambda$. (See Figure \ref{Fig:B1}.)
\end {Prop}

\begin {proof}
For each vertex $V_j \in C$ inside $\Lambda$, the balancing condition
\[
\sum_{k} u_{jk}= 0
\]
holds. Thus
\[
\sum_{j,k} \ddvec{P_0 V_j} \times u_{jk}= 0.
\]
On the left hand side, two weighted primitive vectors of the same edge inside $\Lambda$ are canceled as follows.
\begin {eqnarray*}
\ddvec{P_0 V_j} \times u_{jk}+ \ddvec{P_0 V_{j'}} \times u_{j'k'} & = & \ddvec{P_0 V_j} \times u_{jk}+ \ddvec{P_0 V_{j'}} \times (- u_{jk}) \\
& = & \ddvec{V_{j'} V_j} \times u_{jk} \\
& = & 0.
\end {eqnarray*}
Thus we have the required equation.
\end {proof}

\begin {Def}
The {\em multiplicity} of a vertex $V \in C$ is
\[
\xxMult(V;C)= 2 \cdot \xxarea(T_V).
\]
The {\em intersection multiplicity} of $V \in C_1 \cap C_2$ is
\[
\mu_V= \frac{1}{2} \left( \xxMult(V;C_1 \cup C_2)- \xxMult(V;C_1)- \xxMult(V;C_2) \right).
\]
(If $V$ is not a vertex of $C$, we put $\xxMult(V;C)= 0$.) \\
The formal sum
\[
C_1 \cdot C_2= \sum_{V \in C_1 \cap C_2} \mu_V V
\]
is called the {\em stable intersection} of $C_1$ and $C_2$.
\end {Def}

The stable intersection is characterized as the limit of the transversal intersection (See \cite{First Steps}, Theorem 4.3). If $V \in C_1 \cap C_2$ is a transversal intersection point, this definition is simplified to
\[
\mu_V= |u_E \times u_F|,
\]
where $E \subset C_1, F \subset C_2$ are edges passing through $V$.

\begin {Thm}[Tropical Bezout's Theorem] \label {Thm:A4}
Let $C_1, C_2$ be tropical curves in $\rr^2$. Then the following formula holds.
\[
\deg(C_1 \cdot C_2)= \xxarea( \Delta(C_1)+ \Delta(C_2))- \xxarea( \Delta(C_1))- \xxarea( \Delta(C_2)),
\]
where $\Delta(C_1)+ \Delta(C_2)$ is the Minkowski sum.
\end {Thm}

\begin {proof}
This follows from the above definition and Proposition \ref{Thm:A1}.
\end {proof}

For example, let $c, d \geq 1$, and suppose $\Delta(C_1)= \xxConv \{(0,0), (c,0), (0,c) \}$, $\Delta(C_2)= \xxConv \{(0,0), (d,0), (0,d) \}$ (In algebraic approach, $C_2$ is said to be {\em projective} of degree $d$). Then $\Delta(C_1)+ \Delta(C_2)= \xxConv \{(0,0), (c+d,0), (0,c+d) \}$, and $\deg(C_1 \cdot C_2)= \frac{1}{2} (c+ d)^2- \frac{1}{2} c^2- \frac{1}{2} d^2= cd$. \par

\section {Proof of the surjectivity} \label {Sec:B2}

First we show that $\ph$ is well-defined (Lemma \ref{Prop:A5}, Lemma \ref{Prop:A6}).

\begin {Lem} \label {Prop:A2}
Let $u \in \zz^2$ be a primitive vector. Then for given $\e>0$, there is a primitive vector $v \in \zz^2$ such that
\[
u \times v= 1, \quad |\xxangle(u)- \xxangle(v)|< \e,
\]
where $\xxangle(u)$ denotes the angle of $u$.
\end {Lem}

This lemma is easy.

\begin {Lem} \label {Prop:A3}
Let $E$ be an edge of $C$, and let $P, P', Q, Q'$ be points of E such that $\dvec{PP'}$ $= \dvec{QQ'}$. Then $P'- P \sim Q'- Q$.
\end {Lem}

\begin {proof}
(Figure \ref{Fig:A5}) We may assume that $P', Q$ lie on the interval $[P, Q']$, and that $Q'$ lies in the interior $\xxInt(E)$. Let $v_1, v_2, v_3$ be primitive vectors such that
\begin {equation}
| u_E \times v_i |= 1 \quad (i=1, 2, 3), \label {eq:A3-1}
\end {equation}
\[
\xxangle(u_E)- \e< \xxangle(v_1)< \xxangle(u_E)< \xxangle(v_2)< \xxangle(u_E)+ \e< \xxangle(v_3).
\]
Then we have a triangle, with vertices $P$ and $R_1, R_2 \in \rr^2$, such that
\begin {eqnarray*}
&& \ddvec{PR_i} \mbox{ has direction } v_i \quad (i=1, 2), \\
&& \ddvec{R_1R_2} \mbox{ has direction } v_3, \\
&& Q' \in [R_1, R_2].
\end {eqnarray*}
Take $\e>0$ small enough so that this triangle is disjoint from $C \setminus E$. \par
    Now we take two tropical curves $L, M$ as follows. $L$ consists of one vertex $P$ and three rays $L_0, L_1, L_2 \subset \rr^2$, with $R_1 \in L_1$, $R_2 \in L_2$. $L_0$ is parallel to $E$. $\xxwt(L_1)= \xxwt(L_2)= 1$. (The balancing condition at $P$ follows from (\ref{eq:A3-1}).) $M$ consists of one finite edge $M_0 \subset \rr^2$, four rays $M_1, M_2, M_3, M_4 \subset \rr^2$, and two vertices $V_1, V_2 \in \rr^2$. $M_0$ is parallel to $[R_1, R_2]$, and passes through $Q$. $V_1 \in [P, R_1]$, $V_2 \in [P, R_2]$. For $i=1, 2$, $M_i$ has vertex $V_i$ and passes through $R_i$. $M_3, M_4$ are parallel to $E$. $\xxwt(M_0)= \xxwt(M_1)= \xxwt(M_2)= 1$. \par
    Move $L$ by $\dvec{PP'}$, and denote it by $L'$. Move $M$ by $\dvec{QQ'}$, and denote it by $M'$. Then we have a relation
\[
(C \cdot L- P)- (C \cdot L'- P')= (C \cdot M- Q)- (C \cdot M'- Q').
\]
Thus $P'-P$ is linearly equivalent to $Q'-Q$.
\end {proof}

\begin {figure}[ht]
\begin {center} \input {./picture/bhA5.tex} \end {center}
\caption {}
\label {Fig:A5}
\end {figure}

\begin {Cor} \label {Prop:A4}
Let $E$ be any edge, and suppose that all interior points of $E$ are linearly equivalent. Then all points of $E$ are linearly equivalent.
\end {Cor}

\begin {Lem} \label {Prop:A5}
Let $E$ be a ray of $C$. Then all points of $E$ are linearly equivalent.
\end {Lem}

\begin {proof}
(Figure \ref{Fig:A6}, left) Let $P, Q \in \xxInt(E)$. Take a primitive vector $v$ so that
\[
| u_E \times v |= 1, \quad | \xxangle(u_E)- \xxangle(v) |< \e.
\]
There is a small parallelogram $R_1 R_2 R_3 R_4$ such that
\begin {eqnarray*}
&& \ddvec{R_3 R_1}, \ddvec{R_4 R_2} \mbox{ have direction } v, \\
&& \ddvec{R_2 R_1}, \ddvec{R_4 R_3} \mbox{ have direction } w:= u_E- v, \\
&& P \in [R_1, R_3],~ Q \in [R_2, R_4].
\end {eqnarray*}
Take $\e>0$ small enough so that this parallelogram is disjoint from $C \setminus E$. \par
Let $M_1$ be a tropical curve, consisting of one vertex $R_1$ and three rays $L_0, L_1, L_2 \subset \rr^2$, such that $L_0$ has direction $u_E$, $R_3 \in L_1$, $R_2 \in L_2$, $\xxwt(L_0)= \xxwt(L_1)= 1$. Move $M_1$ by $\dvec{R_1 R_i}$, and denote it by $M_i (i=2, 3, 4)$. Then we have a relation
\[
C \cdot M_1+ C \cdot M_4- Q= C \cdot M_2+ C \cdot M_3- P.
\]
Thus $P \sim Q$.
\end {proof}

\begin {figure}[ht]
\begin {center} \input {./picture/bhA6.tex} \end {center}
\caption {}
\label {Fig:A6}
\end {figure}

\begin {Lem} \label {Prop:A6}
Let $E$ be a tentacle edge of $C$. Then all points of E are linearly equivalent.
\end {Lem}

\begin {proof}
(Figure \ref{Fig:A6}, right) Let $P \in \xxInt(E)$. Since $C \setminus \xxInt(E)$ is disconnected, $E$ is a boundary of two unbounded convex open sets $U_1, U_2$ of $\rr^2 \setminus C$. For $i=1, 2$, let $\theta_i$ be the angle of any unbounded direction of $U_i$. For $\delta, \e>0$, let
\[
W_i= \left\{
\begin {array}[c]{rcl}
P+ a+ w & | & a, w \in \rr^2,~ |a|< \delta,~ P+ a \in E, \\
&& | \xxangle(w)- \theta_i|< \e
\end {array}
\right\}, \]\[
W= W_1 \cup W_2.
\]
Take $\delta, \e>0$ small enough so that $W$ intersects $C \setminus E$ only at points of rays of $C$. \par
Let $v$ be a primitive vector such that $u_E \times v= 1$. Take $w, w' \in \zz^2$ so that
\[
w- w'= v, \]\[
| \xxangle(w)- \theta_1 |< \e, \]\[
| \xxangle(w')- \theta_1 |< \e.
\]
Let $L$ be a tropical curve consisting of three rays $L_0, L_1, L_2 \subset \rr^2$ such that
\[
L_1 \cup L_2 \subset W, \]\[
u_{L_1}= w, \]\[
| \xxangle(L_2)- \theta_2 |< \e, \]\[
L_1 \cap E= P.
\]
Moving $L$ on the direction of $L_0$, the intersection point $P$ changes to other point $Q \in E$, but all other intersection points of $L$ and $C$ are stable except for points of rays of $C$. We have
\[
(u_E \times w) P \sim (u_E \times w) Q.
\]
Similarly,
\[
(u_E \times w') P \sim (u_E \times w') Q.
\]
Thus $P \sim Q$.
\end {proof}

For a proof of Theorem \ref{Prop:A1}, the case of genus $1$ is essential. If $\Lambda := \xxBunch(C) \approx S^1$, then we have a map $\pi \colon \rr \rightarrow \Lambda$ with the following properties. \\
i) $\pi(0)= \oo$. \\
ii) $\pi$ is increasing with respect to the positive direction of $S^1$. \\
iii) If $\pi [a, b]$ ($a, b \in \rr$) is contained in an edge of primitive vector $u$, then
\[
\xxlength( \pi [a, b])=| u|( b- a).
\]
(In other words, $\pi$ is compatible with the lattice length.)

\begin {Lem} \label {Prop:A7}
If $\pi ( a)= P$, $\pi ( a')= P'$, $\pi ( b)= Q$, $\pi ( b')= Q'$, $a'- a= b'- b$, then $P'- P \sim Q'- Q$.
\end {Lem}

\begin {proof}
(Figure \ref{Fig:A7}) We may assume that $a'- a> 0$ is small enough so that $P, P'$ lie on the same edge $E$, and $Q, Q'$ on $F$. We may also assume that $E, F$ are adjacent at a common vertex $R$, and $E \not= F$. From Lemma \ref{Prop:A3}, we may assume that $P, P', Q, Q'$ are interior points of edges, and $[P, Q']$ has rational slope. \par
Let $L$ be the line passing through $P, Q'$. Then
\[
E \cdot L= |u_E \times u_L| P, \]\[
F \cdot L= |u_F \times u_L| Q.
\]
Let $P'', Q''$ be points of $E, F$ such that
\[
\ddvec{P P''}= \frac{1}{|u_E \times u_L|} \ddvec{P P'}, \]\[
\ddvec{Q' Q''}= \frac{1}{|u_F \times u_L|} \ddvec{Q' Q}.
\]
Then
\begin {eqnarray*}
|\ddvec{P P''} \times u_L| & = & \frac{| \ddvec{P P''}|}{| u_E|} | u_E \times u_L| \\
& = & \frac{| \ddvec{P P'}|}{| u_E|} \\
& = & \frac{| \ddvec{Q' Q}|}{| u_F|} \quad ( \mbox{because } a'- a= b'- b) \\
& = & |\ddvec{Q' Q''} \times u_L|,
\end {eqnarray*}
which means that $[P'', Q'']$ is parallel to $L$. Thus
\[
| u_E \times u_L|( P''- P) \sim |u_F \times u_L|( Q'- Q'').
\]
This means $P'- P \sim Q'- Q$, from Lemma \ref{Prop:A3}.
\end {proof}

\begin {figure}[ht]
\begin {center} \input {./picture/bhA7.tex} \end {center}
\caption {}
\label {Fig:A7}
\end {figure}

\begin {Lem}[Interval divisor] \label {Prop:A8}
Let $S, T$ be rays in $\rr^2$ such that $S \cap C= \emptyset$, $T \cap C= \emptyset$, $S \cap T= \emptyset$. Let $P, P' \in S$, $Q, Q' \in T$ be points such that $[P, Q]$ is parallel to $[P', Q']$. Then
\[
C \cdot [P, Q] \sim C \cdot [P', Q'].
\]
\end {Lem}

\begin {proof}
(Figure \ref{Fig:A8}, left) Let $L$ be the tropical curve with vertex $P$, consisting of three rays $L_0, L_1, L_2$ such that
\[
L_1 \subset S, \quad Q \in L_2, \]\[
\xxwt(L_1)= \xxwt(L_2)= 1.
\]
Let $R$ be the point such that
\[
[P, R] \mbox{ is parallel to } L_0, \]\[
[P', R] \mbox{ is parallel to } [P, Q].
\]
Then we have a triangle $P P' R$. We may assume that $|P P'|$ is small enough so that this triangle is disjoint from $C$. Move $L$ by $\dvec{P R}$, and denote it by $L'$. Then we have a relation
\[
C \cdot L- C \cdot [P, Q, \infty)= C \cdot L'- C \cdot [P', Q', \infty).
\]
Thus
\[
C \cdot [P, Q, \infty) \sim C \cdot [P', Q', \infty).
\]
Similarly,
\[
C \cdot [Q, P, \infty) \sim C \cdot [Q', P', \infty), \]\[
C \cdot {\rm Line}(P, Q) \sim C \cdot {\rm Line}(P', Q').
\]
The statement follows.
\end {proof}

\begin {figure}[ht]
\begin {center} \input {./picture/bhA8.tex} \end {center}
\caption {}
\label {Fig:A8}
\end {figure}

\begin {proof}[Proof of the surjectivity of $\ph$.]
(Figure \ref{Fig:A8}, right) Since the image of $\Lambda_1 \cup \cdots \cup \Lambda_g$ in $\xxBunch(C)$ is a bouquet, there are $\oo_1, \ldots, \oo_g \sim \oo$ such that $\oo_i \in \Lambda_i$. Because of convexity, there are connected disjoint $g$ cones $U_1, \ldots, U_g \subset \rr^2$ with center $\oo_1, \ldots, \oo_g$ such that $\Lambda_i \subset U_i$. Similarly to the case of genus $1$, we have a map $\pi_i \colon \rr \rightarrow \Lambda_i$ for each $i$. Lemma \ref{Prop:A7} is proved for $\pi_i$ similarly, only changing $L$ to an interval divisor $L \cap U_i$ (Lemma \ref{Prop:A8}). Thus we have a homomorphism of abelian groups
\[
\tilde{ \ph} \colon \rr^g \underset{ \pi}{ \rmap} \Lambda_1 \times \cdots \times \Lambda_g \underset{ \ph}{ \rmap} \xxJac(C).
\]
Since $\xxJac(C)$ is generated by $\{ P_i- \oo| P_i \in \Lambda_i, 1 \leq i \leq g \}$, $\tilde{ \ph}$ is surjective.
\end {proof}

\section {Parameter space of tropical plane curves}

Let $L$ be a tropical curve with Newton complex $\nnn$. Let $V_0= (b_1, b_2)$ be a fixed vertex of $L$. Let $E_1, \ldots, E_l$ be all finite edges of $L$. Let $a_i$ be the lattice length of $E_i$. Then all tropical curves with Newton complex $\nnn$ are parametrized by $a_1, \ldots, a_l> 0$ and $b_1, b_2 \in \rr$. Let $\xxTN( \nnn, \rr^2) \subset \rr^{l+2}$ be the parameter space.

\begin {Prop} \label {Prop:C1}
$\xxTN( \nnn, \rr^2)$ is connected.
\end {Prop}

\begin {proof}
Let $\Gamma_1, \ldots, \Gamma_g$ be all convex cycles of $L$. Let $E_{i(j,1)}, \ldots, E_{i(j,s_j)}$ be all edges of $\Gamma_j$. Let $u_{j,k}$ be the primitive vector of $E_{i(j,k)}$ of positive direction of $S^1$. Then the equation
\begin {equation}
a_{i(j,1)} u_{j,1}+ \cdots + a_{i(j,s_j)} u_{j,s_j}= 0 \label {eq:C1-1}
\end {equation}
is satisfied for any $L \in \xxTN( \nnn, \rr^2)$. Let $H_j \subset \rr^{l+2}$ be the linear subspace defined by equation (\ref{eq:C1-1}). Then
\[
\xxTN( \nnn, \rr^2)= \{ (a_1, \ldots, a_l, b_1, b_2) | a_1, \ldots, a_l> 0 \} \cap (H_1 \cap \cdots \cap H_g).
\]
Thus $\xxTN( \nnn, \rr^2)$ is a relatively open convex cone in $\rr^{l+2}$, which is connected.
\end {proof}

\begin {Def}
A tropical curve $L'$ is a {\em degeneration} of $L$ if $\Delta(L')= \Delta(L)$ and $\xxNewt(L') \subset \xxNewt(L)$.
\end {Def}

The set of all degenerations of $L$ is parametrized by $\overline{\xxTN( \nnn, \rr^2)}$. If a Newton polygon $\Delta$ is fixed, all tropical curves have a common degeneration (that is, a tropical curve consisting of one vertex). All the Newton polygons is countable. Therefore, the space $\xxTN( \rr^2)$ of all tropical curves is a disjoint union of countable closed cones in affine spaces.

\begin {Cor} \label {Prop:D1}
Tropical curves $L, L' \in \xxTN( \rr^2)$ lie in the same connected component if and only if $\Delta(L)= \Delta(L')$.
\end {Cor}

\section {Proof of the injectivity}

Let $\pi \colon \rr \rightarrow \Lambda$ be the map defined in section \ref{Sec:B2}. $\Lambda$ is considered as a residue group of $\rr$. Let $E_1, \ldots, E_N$ be all edges of $\Lambda$ ordered by the positive direction of $S^1$. Let $\lambda \colon C \rightarrow \Lambda$ be the canonical surjection. For a tropical curve $L \in \xxTN( \rr^2)$, we define $\sigma(L) \in \Lambda$ as follows.
\[
C \cdot L= P_1+ \cdots+ P_r, \]\[
\sigma(L)= \lambda(P_1)+ \cdots + \lambda(P_r).
\]

\begin {Lem} \label {Prop:C2}
$\sigma \colon \xxTN( \rr^2) \rightarrow \Lambda$ is locally constant.
\end {Lem}

\begin {proof}
$\sigma$ is continuous by definition of the stable intersection. Let $\{ L_t | 0 \leq t \leq 1 \}$ be a continuous family of tropical curves with Newton complex $\nnn$ such that $C$ intersects $L_t$ transversely for any $t$. Then there are points $P_{ijt} \in E_i$, edges $L_{ijt} \subset L_t$, and vectors $u_{ij} \in \rr^2$ such that \\
i) $E_i \cdot L_t= \sum_{j} P_{ijt}$, \\
ii) $E_i \cap L_{ijt}= P_{ijt}$, \\
iii) $u_{ij}$ is the weighted primitive vector of $L_{ijt}$ starting at the vertex inside $\Lambda$, divided by $\mu_{P_{ijt}}$. \\
For $L_0$ and $L_1$, we have the moment condition inside $\Lambda$:
\[
\sum_{i,j} \xxmoment(u_{ij}, P_{ij0})= 0, \]\[
\sum_{i,j} \xxmoment(u_{ij}, P_{ij1})= 0.
\]
From these,
\[
\sum_{i,j} \left( \ddvec{P_{ij0} P _{ij1}} \times u_{ij} \right)= 0.
\]
Since $u_{E_i} \times u_{ij}= -1$, this means
\[
\sum_{i,j} \left( \lambda(P_{ij0})- \lambda(P_{ij1}) \right)= 0.
\]
Thus $\sigma(L_0)= \sigma(L_1)$.
\end {proof}

\begin {figure}[ht]
\begin {center} \input {./picture/bhC3.tex} \end {center}
\caption {}
\label {Fig:C3}
\end {figure}

\begin {proof}[Proof of the injectivity of $\ph$.]
(Figure \ref{Fig:C3}) Suppose
\[
(P_1+ \cdots + P_g)- (Q_1+ \cdots + Q_g)= C \cdot L- C \cdot L', \]\[
P_i, Q_i \in \Lambda_i \setminus \xxVer( \Lambda_i), \]\[
\Delta(L)= \Delta(L').
\]
Let $\tilde{C}$ be the tropical curve consisting of $\Lambda_1$ and $N$ rays $F_1, \ldots, F_N$. Then $\deg(F_i \cdot L)= \deg(F_i \cdot L')$ because of the tropical Bezout's theorem. From Corollary \ref{Prop:D1} and Lemma \ref{Prop:C2}, we have $\sigma(L)= \sigma(L')$. Thus $P_1= Q_1$.
\end {proof}

\end {document}

%% file: picture/bhA1.tex
\begin {picture}(0, 0)
\includegraphics {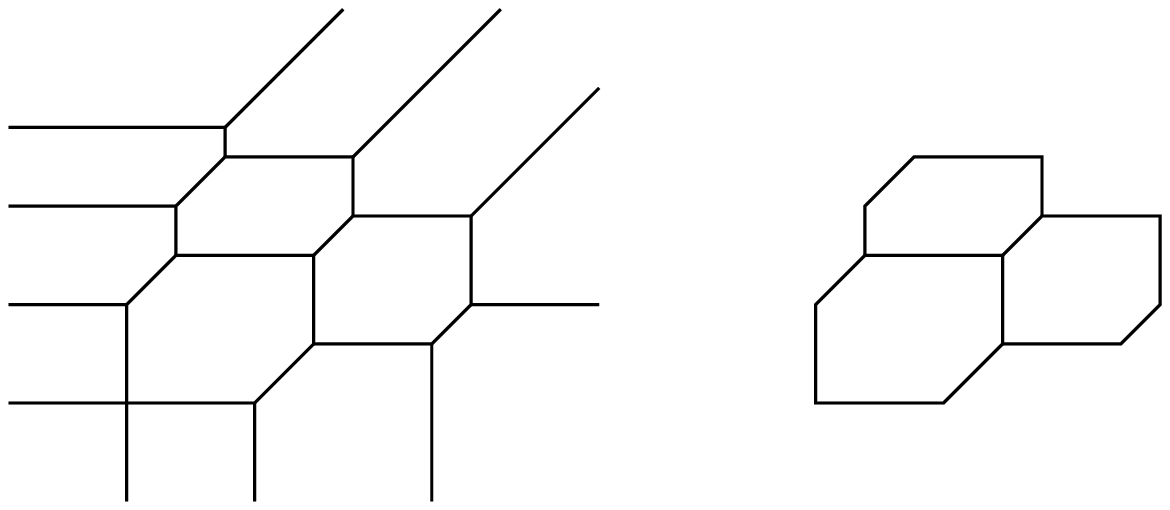}
\end {picture}
\begin {picture}(350, 150)
\unitlength= 1 mm

\put (28, -5) {$C$}
\put (85, -5) {$\xxBunch(C)$}

\end {picture}

%% file: picture/bhA3.tex
\begin {picture}(0, 0)
\includegraphics {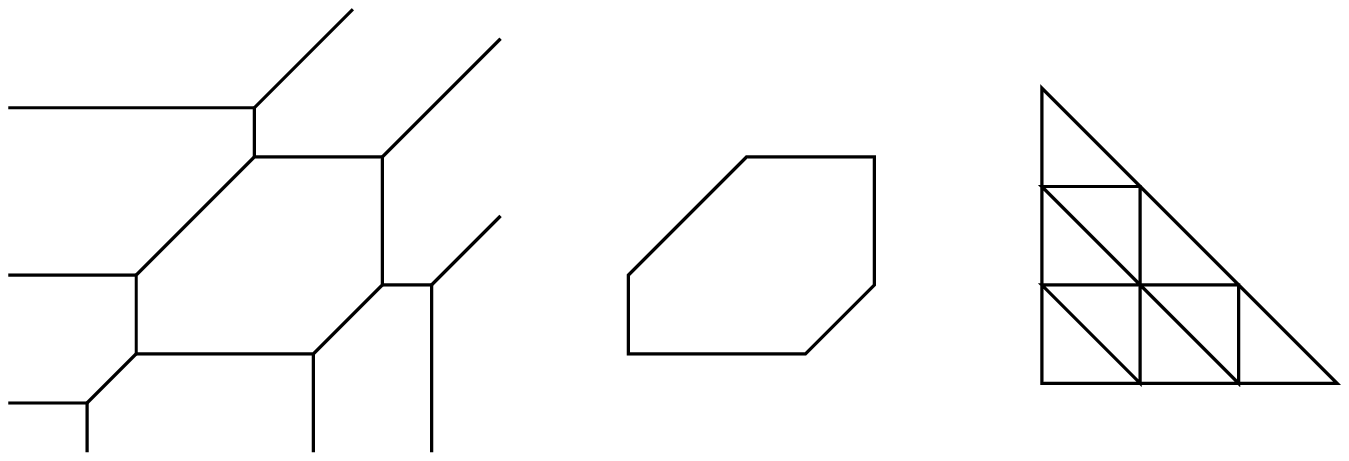}
\end {picture}
\begin {picture}(350, 150)
\unitlength= 1 mm

\put (20, 0) {$C$}
\put (62, 0) {$\xxBunch(C)$}
\put (106, 0) {$\xxNewt(C)$}

\end {picture}

%% file: picture/bhA4.tex
\begin {picture}(0, 0)
\includegraphics {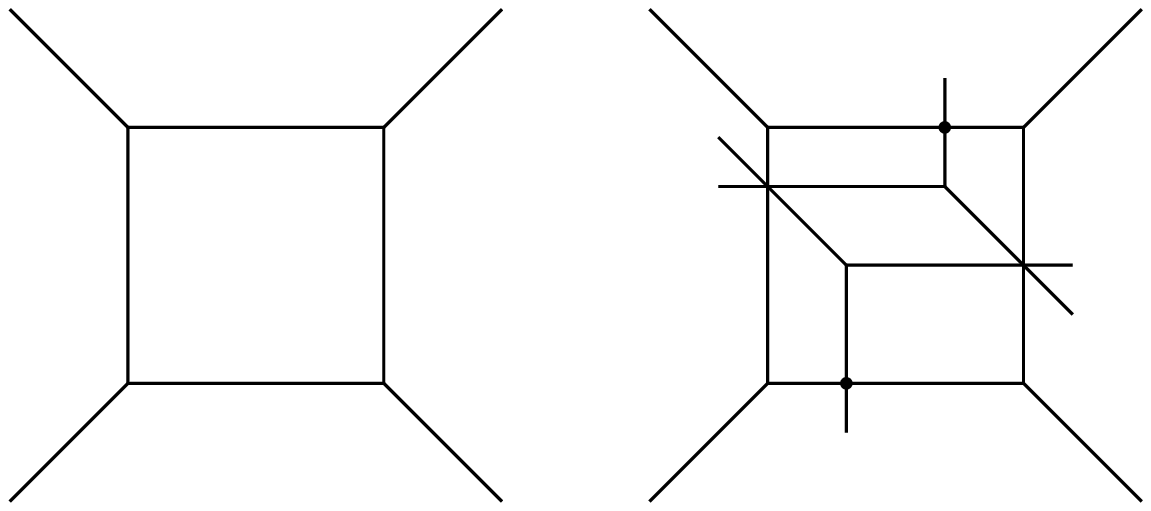}
\end {picture}
\begin {picture}(350, 150)
\unitlength= 1 mm

\put (25, -5) {$C$}
\put (90, 8) {$P$}
\put (100, 40) {$Q$}
\put (86, 2) {$L$}
\put (96, 46) {$M$}

\end {picture}

%% file: picture/bhA2.tex
\begin {picture}(0, 0)
\includegraphics {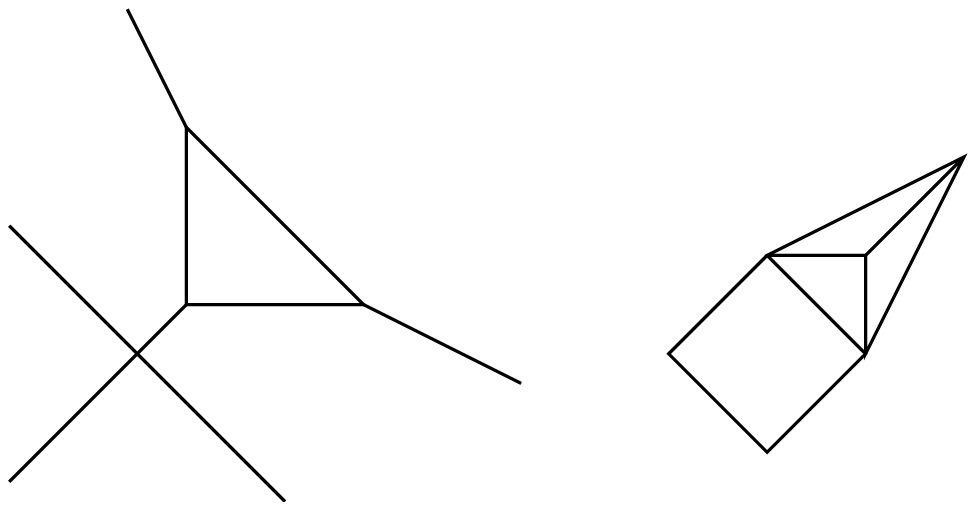}
\end {picture}
\begin {picture}(350, 150)
\unitlength= 1 mm

\put (28, -5) {$C$}
\put (82, -5) {$\xxNewt(C)$}

\end {picture}

%% file: picture/bhB1.tex
\begin {picture}(0, 0)
\includegraphics {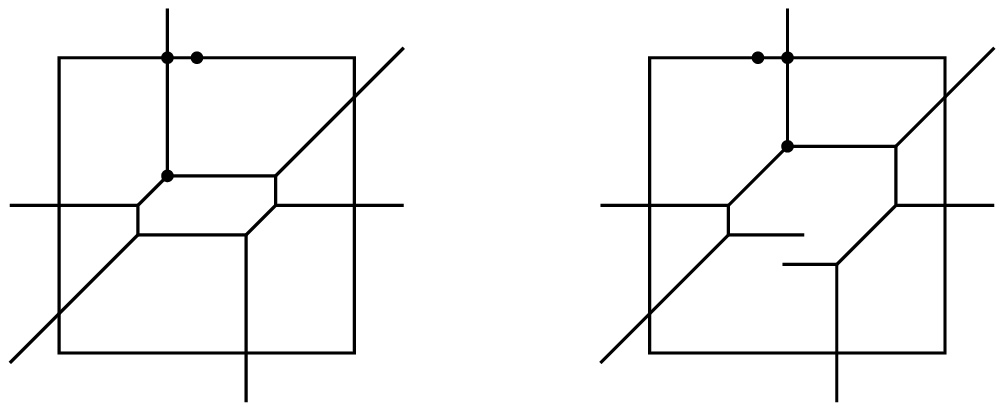}
\end {picture}
\begin {picture}(350, 150)
\unitlength= 1 mm

\put (20, 41) {$P$}
\put (80, 41) {$P$}
\put (29, 41) {$P'$}
\put (89, 41) {$P'$}
\put (23, 47) {$E_1$}
\put (86, 47) {$E'_1$}
\put (17, 29) {$V_{E_1}$}
\put (80, 32) {$V_{E'_1}$}

\end {picture}

%% file: picture/bhA5.tex
\begin {picture}(0, 0)
\includegraphics {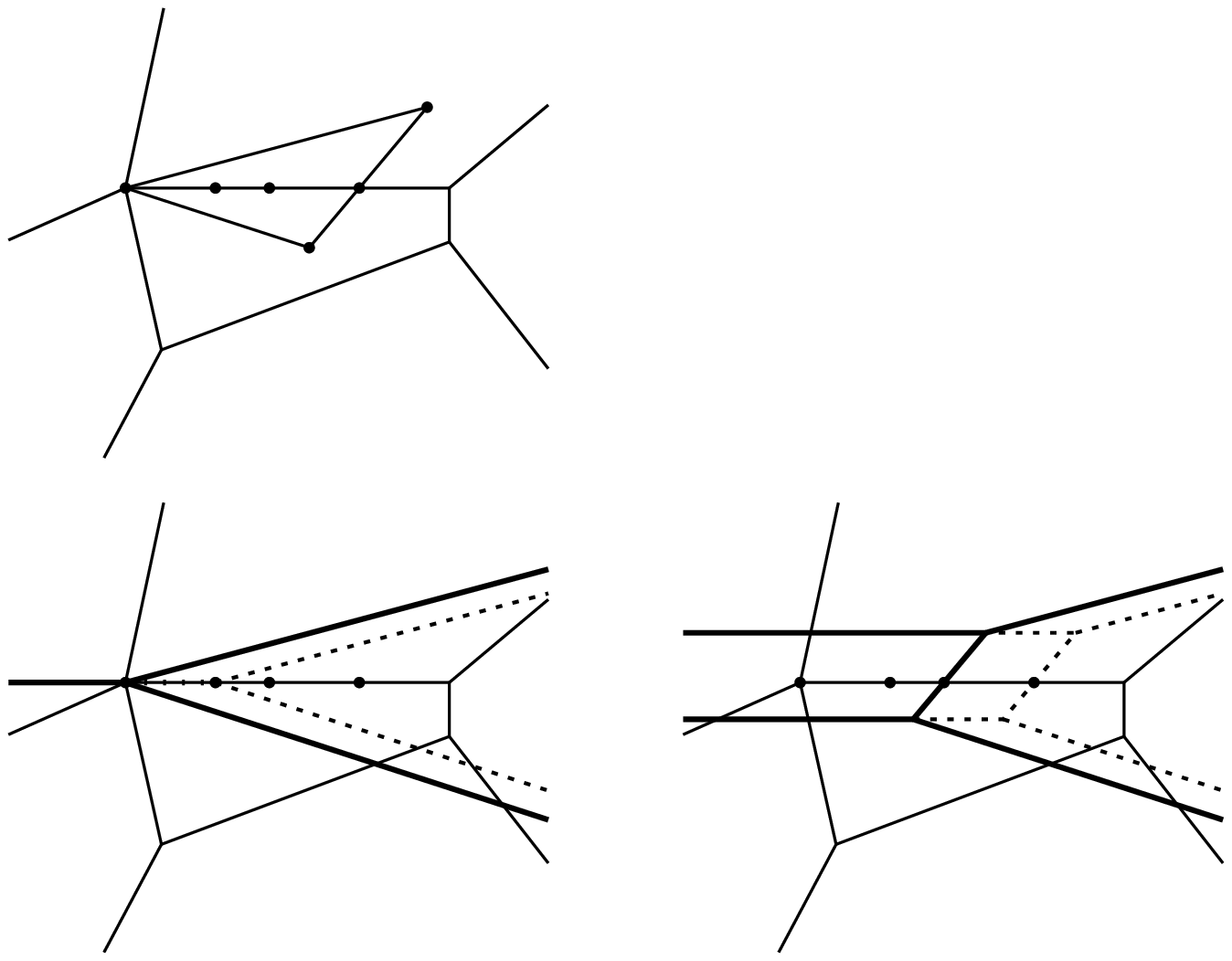}
\end {picture}
\begin {picture}(350, 300)
\unitlength= 1 mm

\put (3, 86) {$P$}
\put (12, 88) {$P'$}
\put (26, 86) {$Q$}
\put (36, 86) {$Q'$}
\put (21, 75) {$R_1$}
\put (38, 96) {$R_2$}

\put (3, 31) {$P$}
\put (12, 33) {$P'$}
\put (-5, 31) {$L$}

\put (101, 31) {$Q$}
\put (111, 31) {$Q'$}
\put (70, 36) {$M$}

\end {picture}

%% file: picture/bhA6.tex
\begin {picture}(0, 0)
\includegraphics {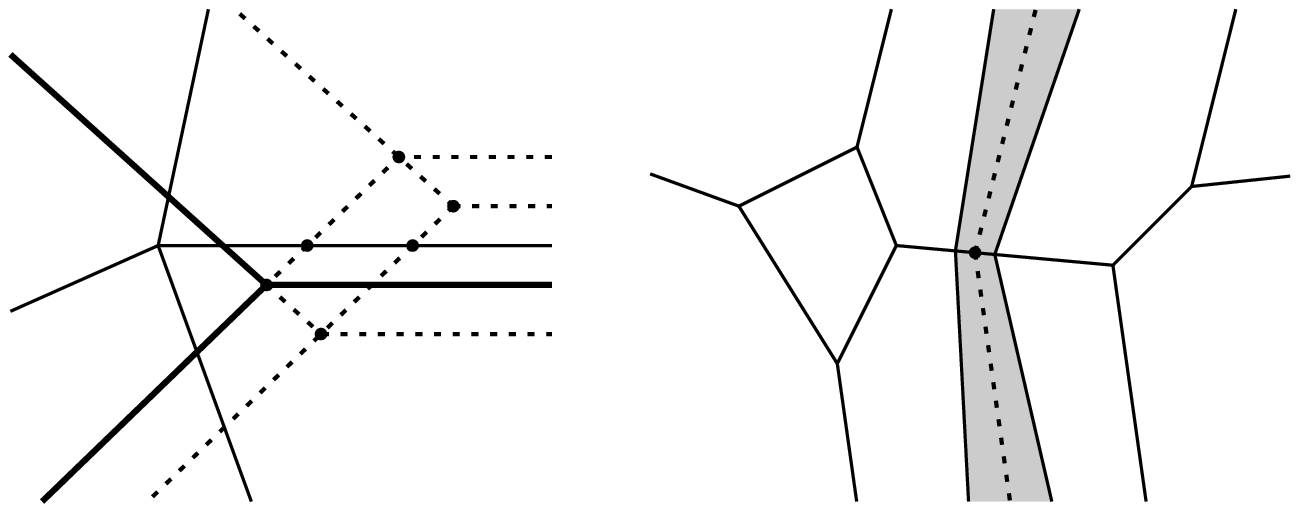}
\end {picture}
\begin {picture}(350, 150)
\unitlength= 1 mm

\put (44, 31) {$R_1$}
\put (39, 36) {$R_2$}
\put (31, 12) {$R_3$}
\put (23, 16) {$R_4$}
\put (3, 43) {$M_4$}
\put (55, 25) {$E$}
\put (39, 22) {$P$}
\put (25, 28) {$Q$}

\put (100, -5) {$\theta_2$}
\put (102, 52) {$\theta_1$}
\put (101, 26) {$P$}

\end {picture}

%% file: picture/bhA7.tex
\begin {picture}(0, 0)
\includegraphics {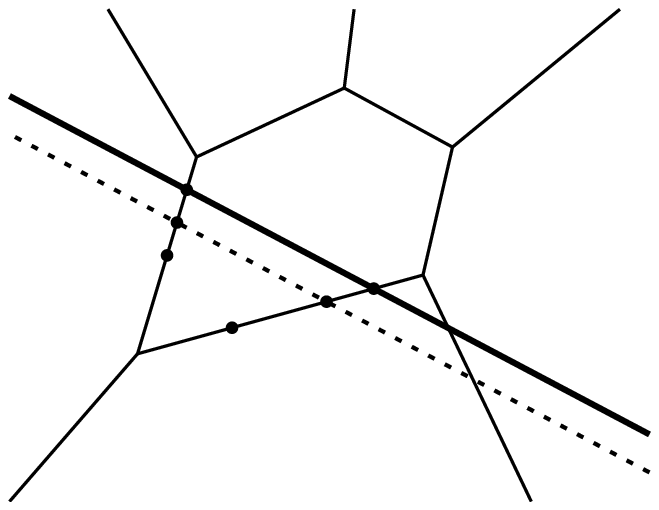}
\end {picture}
\begin {picture}(350, 150)
\unitlength= 1 mm

\put (90, 8) {$L$}
\put (44, 32) {$P$}
\put (32, 27) {$P''$}
\put (34, 22) {$P'$}

\put (44, 13) {$Q$}
\put (53, 16) {$Q''$}
\put (59, 25) {$Q'$}

\end {picture}

%% file: picture/bhA8.tex
\begin {picture}(0, 0)
\includegraphics {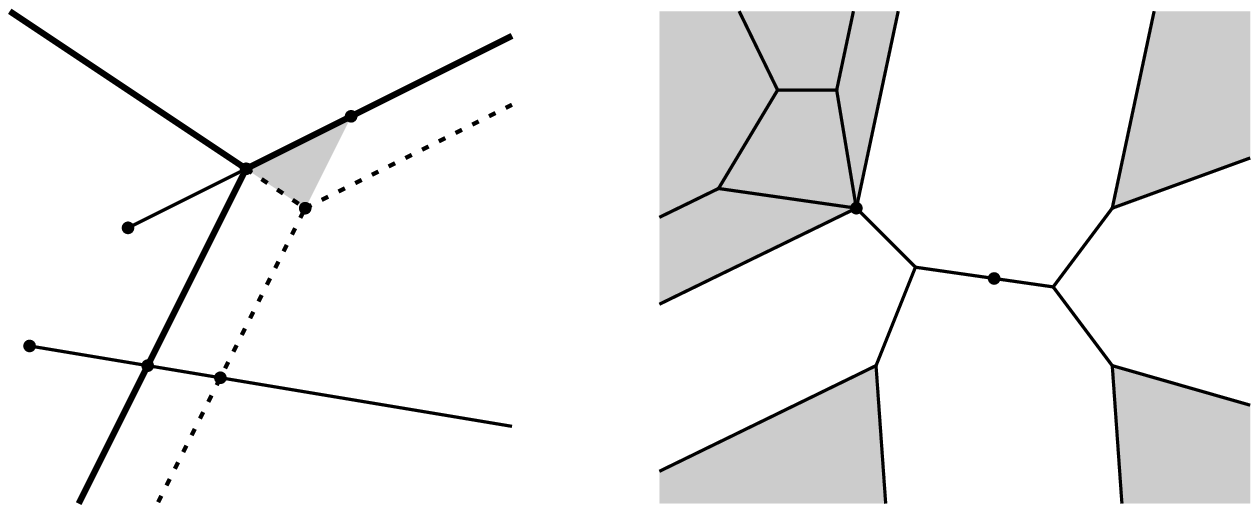}
\end {picture}
\begin {picture}(350, 150)
\unitlength= 1 mm

\put (26, 37) {$P$}
\put (35, 40) {$P'$}
\put (33, 25) {$R$}
\put (11, 10) {$Q$}
\put (25, 8) {$Q'$}
\put (55, 46) {$S$}
\put (55, 5) {$T$}
\put (4, 42) {$L$}

\put (102, 24) {$\oo$}
\put (91, 29) {$\oo_i$}

\end {picture}

%% file: picture/bhC3.tex
\begin {picture}(0, 0)
\includegraphics {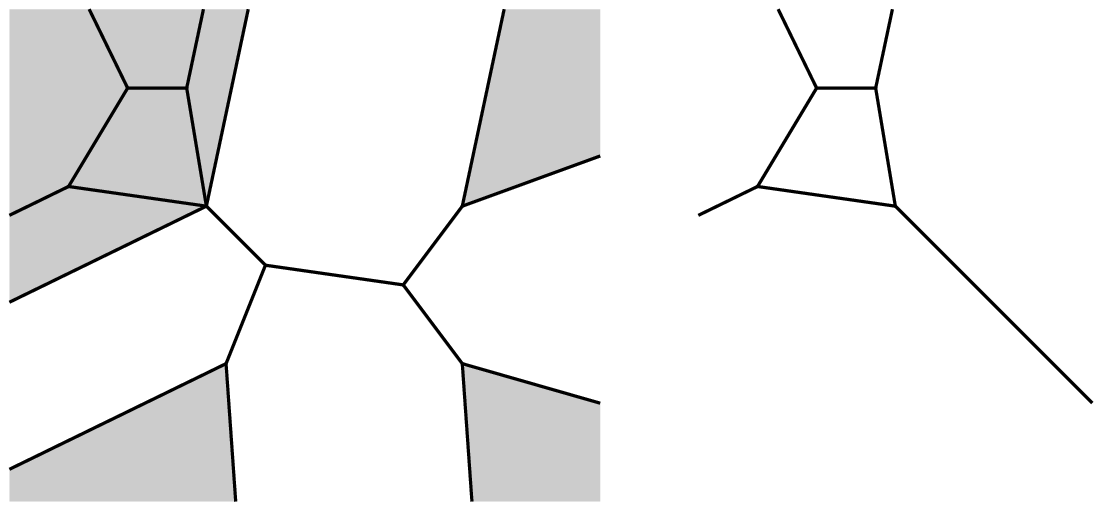}
\end {picture}
\begin {picture}(350, 150)
\unitlength= 1 mm

\put (30, -5) {$C$}
\put (90, -5) {$\tilde{C}$}
\put (108, 12) {$F_i$}

\end {picture}